\documentclass[11pt]{amsart}

\usepackage{graphicx,pinlabel}

\newtheorem{theorem}{Theorem}[section]
\newtheorem{proposition}[theorem]{Proposition}

\newtheorem{corollary}[theorem]{Corollary}

\theoremstyle{definition}
\newtheorem{definition}[theorem]{Definition}

\numberwithin{equation}{section}

\newcommand{\Q}{\operatorname{{\mathbb Q}}}
\newcommand{\R}{\operatorname{{\mathbb R}}}
\newcommand{\Z}{\operatorname{{\mathbb Z}}}

\newcommand{\G}{\operatorname{{\mathcal G}}}

\begin{document}

\title[Cabling sequences of tunnels of torus knots]
{Cabling sequences of tunnels of torus knots}

\author{Sangbum Cho}
\address{Department of Mathematics\\
University of California at Riverside\\
Riverside, California 92521\\
USA} 
\email{scho@math.ucr.edu}

\author{Darryl McCullough}
\address{Department of Mathematics\\
University of Oklahoma\\
Norman, Oklahoma 73019\\
USA} 
\email{dmccullough@math.ou.edu}
\urladdr{www.math.ou.edu/$_{\widetilde{\phantom{n}}}$dmccullough/}
\thanks{The second author was supported in part by NSF grant DMS-0802424}


\date{\today}


\begin{abstract} 
In previous work, we developed a theory of tunnels of tunnel number $1$
knots in $S^3$. It yields a parameterization in which each tunnel is
described uniquely by a finite sequence of rational parameters and a finite
sequence of $0$'s and $1$'s, that together encode a procedure for
constructing the knot and tunnel. In this paper we calculate these
invariants for all tunnels of torus knots.
\end{abstract}

\maketitle

\section*{Introduction}
\label{sec:intro}

In previous work~\cite{CM}, we developed a theory of tunnels of tunnel
number~$1$ knots in $S^3$. It shows that every tunnel can be obtained from
the unique tunnel of the trivial knot by a uniquely determined sequence of
``cabling constructions''. A cabling construction is determined by a
rational parameter, called its ``slope,'' so this leads to a
parameterization of all tunnels of all tunnel number~$1$ knots by sequences
of rational numbers and ``binary'' invariants. Various applications of the
theory are given in~\cite{CM}, \cite{CMsteps} and~\cite{CMdepth}, as well
as other work in preparation.

Naturally, it is of interest to calculate these invariants for known
examples of tunnels. In~\cite{CM}, they are calculated for all tunnels of
$2$-bridge knots, and in the present paper we obtain them for all tunnels
of torus knots. Tunnels of torus knots are a key example in our study of
the ``depth'' invariant in~\cite{CMdepth}. Also, torus knots are special in
that their complements have zero (Gromov) volume, so they should be
critical to understanding how hyperbolic volumes of complements of tunnel
number~$1$ knots are related to the sequences of slope and binary
invariants of their tunnels.

In the next section, we will give the main results.
Sections~\ref{sec:tunnels}, \ref{sec:slopedisks}, and \ref{sec:cabling}
provide a concise review of the parts of the theory from~\cite{CM} that
will be needed in this paper. The main results are proven in
Section~\ref{sec:regular} for the middle tunnels and
Section~\ref{sec:semisimple} for the upper and lower tunnels.

The calculations in this paper enable us to recover the classification of
torus knot tunnels given by M. Boileau, M. Rost, and
H. Zieschang~\cite{B-R-Z} and Y. Moriah~\cite{Moriah}, although not their
result that these are all the tunnels. We give this application in
Section~\ref{sec:applications} below.

All of our algorithms to find the invariants are straightforward to
implement computationally, and we have done this in software available at
\cite{slopes}. Sample computations are given in
Section~\ref{sec:statement}.

In work in progress, we are developing a general method for computing these
invariants for all $(1,1)$-tunnels. In particular, this will recover the
calculations for tunnels of $2$-bridge knots, given in~\cite{CM}, and for
some of the tunnels of torus knots that we give here (the upper and lower
tunnels, but not the middle tunnels). Still, we think it is worthwhile to
give the method of this paper, which is more direct and more easily
visualized.

We are grateful to the referee for a prompt and careful reading of the
original manuscript.

\section{The main results}
\label{sec:statement}

To set notation, consider a (nontrivial) $(p,q)$ torus knot $K_{p,q}$,
contained in a standard torus $T$ bounding a solid torus $V\subset
S^3$. In $\pi_1(V)$, $K_{p,q}$ represents $p$ times a generator. The
complementary torus $\overline{S^3-V}$ will be denoted by~$W$.

The tunnels of torus knots were classified by M. Boileau, M. Rost, and
H. Zieschang \cite{B-R-Z} and Y. Moriah~\cite{Moriah}.  The \textit{middle
tunnel} of $K_{p,q}$ is represented by an arc in $T$ that meets $K_{p,q}$
only in its endpoints. The \textit{upper tunnel} of $K_{p,q}$ is
represented by an arc $\alpha$ properly imbedded in $W$, such that the
circle which is the union of $\alpha$ with one of the two arcs of $K_{p,q}$
with endpoints equal to the endpoints of $\alpha$ is a deformation retract
of $W$. The \textit{lower tunnel} is like the upper tunnel, but
interchanging the roles of $V$ and~$W$. In certain cases, some of these
tunnels are equivalent, as we will detail in
Section~\ref{sec:applications}.

To state our results for the middle tunnels, assume for now that
$p,q>1$. Since $K_{p,q}$ and $K_{q,p}$ are equivalent by an
orientation-preserving homeomorphism of $S^3$ taking middle tunnel to
middle tunnel, we may also assume that $p>q$. Put
$U=\begin{pmatrix}1&1\\0&1\end{pmatrix}$ and
$L=\begin{pmatrix}1&0\\1&1\end{pmatrix}$.
\begin{theorem} Let $p$ and $q$ be relatively prime integers with $p>q\geq
2$. Write $p/q$ as a continued fraction $[n_1,n_2,\ldots, n_k]$ with all
$n_j$ positive and $n_k\neq 1$. Let
\[ A_i = \begin{cases}%
L & -n_1\leq i\leq -1\\
U & 0\leq i\leq n_2-1\\
L & n_2\leq i\leq n_2+n_3-1\\
U & n_2+n_3\leq i\leq n_2+n_3+n_4-1\\
& \cdots \\
L& k\text{\ odd and\ }n_2+n_3+\cdots +n_{k-1}\leq i\leq n_2+n_3+\cdots
+ n_k-1\\
U& k\text{\ even and\ }n_2+n_3+\cdots +n_{k-1}\leq i\leq n_2+n_3+\cdots
+ n_k-1\ .
\end{cases}
\]
Put $N=n_2+n_3+\cdots + n_k-2$, and for $0\leq t\leq N$ put
\[\begin{pmatrix} a_t & b_t\\ c_t & d_t\end{pmatrix}=\prod_{i=t}^{-n_1}A_i
\ ,\]
where the subscripts in the product occur in descending order. Then:
\begin{enumerate}
\item[(i)] The middle tunnel of $K_{p/q}$ is
produced by $N+1$ cabling constructions whose slopes $m_0$, $m_1,\ldots\,$,
$m_N$ are
\[ \left[\frac{1}{2n_1+1}\right], \;a_1d_1+b_1c_1,\; a_2d_2+b_2c_2,\; \ldots,\;
a_Nd_N+b_Nc_N\ .\]
\item[(ii)] For each $t$, the cabling corresponding to the slope invariant
$m_t$ produces the $(a_t+c_t,b_t+d_t)$ torus knot; in particular, the first
cabling produces the $(2n_1+1,2)$ torus knot.
\item[(iii)] The binary invariants of the cabling sequence of this tunnel,
for $2\leq t\leq N$, are given by $s_t=1$ if $A_t\neq A_{t-1}$ and $s_t=0$
otherwise.
\end{enumerate}
\label{thm:middle_tunnels_slopes}
\end{theorem}
\noindent If $pq<0$, then $K_{p,q}$ is equivalent to $K_{|p|,|q|}$ by an
orientation-reversing homeomorphism taking the middle tunnel to the middle
tunnel, so the cabling slopes for the middle tunnel of $K_{p,q}$ are just
the negatives of those of $K_{|p|,|q|}$ given in
Theorem~\ref{thm:middle_tunnels_slopes}, while the binary invariants are
unchanged.

It is not difficult to implement this calculation computationally, and we
have made a script for this available~\cite{slopes}. For $K_{41,29}$, we find
\smallskip

\noindent \texttt{TorusKnots$>$ middleSlopes(41, 29)}\\
\texttt{[ 1/3 ], 5, 17, 29, 99, 169, 577}
\smallskip

\noindent and for $K_{181,-48}$
\smallskip

\noindent \texttt{TorusKnots$>$ middleSlopes(181, -48)}\\
\texttt{[ 6/7 ], -15, -23, -31, -151, -271, -883, -2157, -3431}
\smallskip

\noindent The torus knots that are the intermediate knots in the cabling
sequence are found by
\smallskip

\noindent \texttt{TorusKnots> intermediates( 41, 29 )}\\
\texttt{(3,2), (4,3), (7,5), (10,7), (17,12), (24,17), (41,29)}
\smallskip

\noindent and the binary invariants by
\smallskip

\noindent \texttt{TorusKnots> binaries(41, 29)}\\
\texttt{[1, 0, 1, 0, 1]}
\smallskip

Now we consider the upper and lower tunnels. Since these are semisimple
tunnels, their binary invariants $s_i$ are all $0$ (see
Section~\ref{sec:cabling}). The cabling slopes are given as follows.
\begin{theorem}
Let $p$ and $q$ be relatively prime integers, both greater than $1$. 
For integers $k$ with $1\leq k\leq q$, define integers $p_k$ by 
\[p_k=\lceil kp/q\rceil = \min\{ j\;|\; jq/p \geq k\}\ ,\] 
and let $k_0= \min\{ k\;|\; p_k > 1\}$. Then
the upper tunnel of $K_{p,q}$ is produced by $q-k_0$ cabling operations, whose
slopes are
\[[1/(2p_{k_0}-1)],\; 2p_{k_0+1}-1,\;\ldots\,, \;2p_{q-1}-1\ .\]\par
\label{thm:cabling_sequence}
\end{theorem}
\noindent As before, when $pq<0$ the slopes are just the negatives of those
given in Theorem~\ref{thm:cabling_sequence} for $K_{|p|,|q|}$. The lower
tunnel of $K_{p,q}$ is equivalent to the upper tunnel of $K_{q,p}$, so
Theorem~\ref{thm:cabling_sequence} also finds the slope sequences of all
lower tunnels.

Again, this algorithm is easily scripted and is available
at~\cite{slopes}. Sample calculations are:\par
\smallskip

\noindent 
\texttt{TorusKnots$>$ upperSlopes( 18, 7 )}\\
\texttt{[ 1/ 5 ], 11, 15, 21, 25, 31}
\smallskip

\noindent 
\texttt{TorusKnots$>$ upperSlopes( 7, 18 )}\\
\texttt{[ 1/ 3 ], 3, 3, 5, 5, 7, 7, 7, 9, 9, 11, 11, 11, 13, 13}
\smallskip

\noindent 
\texttt{TorusKnots$>$ lowerSlopes( 18, 7 )}\\
\texttt{[ 1/ 3 ], 3, 3, 5, 5, 7, 7, 7, 9, 9, 11, 11, 11, 13, 13}
\smallskip

Theorems~\ref{thm:middle_tunnels_slopes} and~\ref{thm:cabling_sequence}
show immediately the following integrality result:
\begin{corollary}
Let $\tau$ be a tunnel of a torus knot. Then the first slope invariant
$m_0$ of $\tau$ is of the form $[1/n]$ for some odd integer $n$, and all
other slopes are odd integers.\par
\end{corollary}
\noindent For the middle tunnels, the integrality of the slope invariants
$m_i$ for $i\geq 1$ follows from the work of Scharlemann and
Thompson~\cite{Scharlemann-Thompson} (which inspired our work
in~\cite{CM}). For as shown in~\cite[Section 14]{CM}, their invariant
$\rho(\tau)$ is our final (or ``principal'') slope invariant $m_N$ reduced
modulo $2$ (that is, viewed as an element of $\Q/2\Z$). Scharlemann and
Thompson computed that the $\rho$-invariants of the middle tunnels are $1$,
so it follows that $m_N$ must be an odd integer. As our construction in
Section~\ref{sec:regular} will show, the intermediate slope invariants
$m_i$ are principal slope invariants for middle tunnels of other torus
knots, so they too must be odd integers.

\section{Tunnels as disks}
\label{sec:tunnels}

This section gives a brief overview of the theory in~\cite{CM}. Fix a
standard unknotted handlebody $H$ in $S^3$.  Regard a tunnel of $K$ as a
$1$-handle attached to a neighborhood of $K$ to obtain an unknotted
genus-$2$ handlebody. Moving this handlebody to $H$, a cocore disk for the
$1$-handle moves to a nonseparating disk in $H$.  The indeterminacy due
to the choice of isotopy is exactly the Goeritz group~$\G$, studied in
\cite{Akbas,Cho,ScharlemannTree}. Consequently, the collection of all
tunnels of all tunnel number~$1$ knots, up to orientation-preserving
homeomorphism, corresponds to the orbits of nonseparating disks in $H$
under the action of~$\G$.  From~\cite{Akbas,Cho,ScharlemannTree}, the
action can be understood and the equivalence classes, i.e.~the tunnels,
arranged in a treelike structure which encodes much of the topological
structure of tunnel number~$1$ knots and their tunnels.

When a nonseparating disk $\tau \subset H$ is regarded as a tunnel, the
corresponding knot is a core circle of the solid torus that results from
cutting $H$ along $\tau$. This knot is denoted by~$K_\tau$.

A disk $\tau$ in $H$ is called \textit{primitive} if there is a disk
$\tau'$ in $\overline{S^3-H}$ such that $\partial \tau$ and $\partial
\tau'$ cross in one point in $\partial H$. Equivalently, $K_\tau$ is the
trivial knot in~$S^3$. All primitive disks are equivalent under the action
of~$\G$. This equivalence class is the unique tunnel of the trivial knot.

A \textit{primitive pair} is an isotopy class of two disjoint nonisotopic
primitive disks in $H$. A \textit{primitive triple} is defined similarly.

\section{Slope disks and cabling arcs}
\label{sec:slopedisks}

\begin{figure}
\labellist
\small \hair 2pt
\pinlabel $\lambda$ [B] at -18 136
\pinlabel $\rho$ [B] at 595 136
\endlabellist
\begin{center}
\includegraphics[width=65 ex]{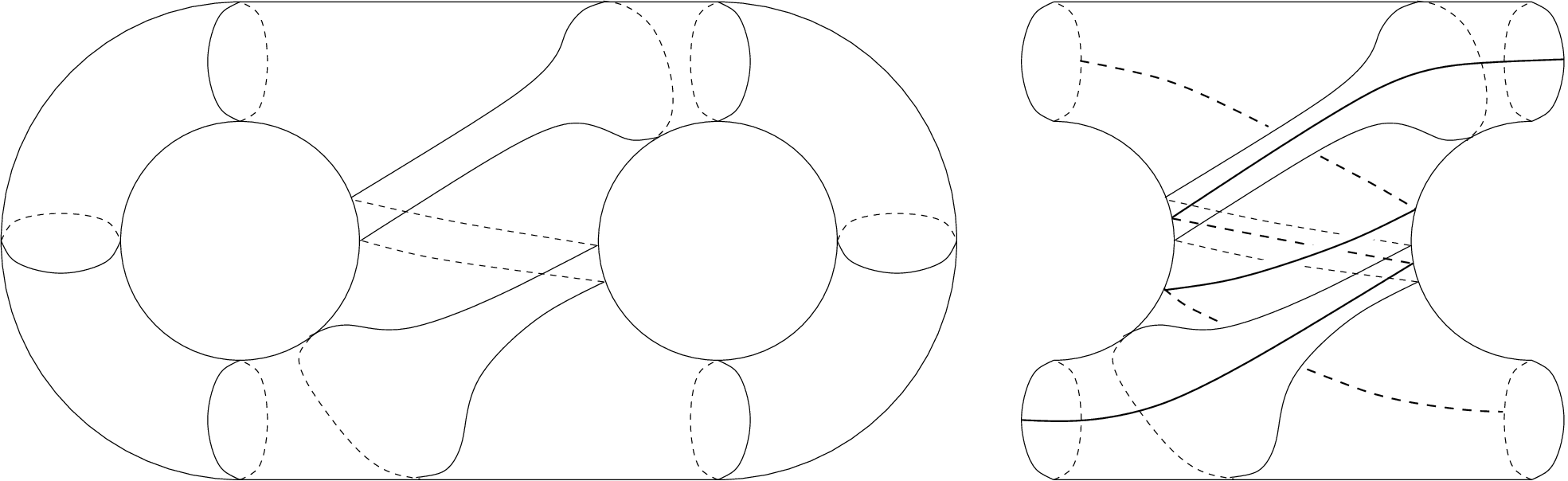}
\caption{A slope disk of $\{\lambda,\rho\}$, and a pair of its cabling arcs
contained in $B$.}
\label{fig:slopedisks}
\end{center}
\end{figure}

This section gives the definitions needed for computing slope invariants.
Fix a pair of nonseparating disks $\lambda$ and $\rho$ (for ``left'' and
``right'') in the standard unknotted handlebody $H$ in $S^3$, as shown
abstractly in Figure~\ref{fig:slopedisks}. The pair $\{\lambda, \rho\}$ is
arbitrary, so in the true picture in $H$ in $S^3$, they will typically look
a great deal more complicated than the pair shown in
Figure~\ref{fig:slopedisks}. Let $N$ be a regular neighborhood of
$\lambda\cup \rho$ and let $B$ be the closure of $H-N$. The frontier of $B$
in $H$ consists of four disks which appear vertical in
Figure~\ref{fig:slopedisks}. Denote this frontier by $F$, and let $\Sigma$
be $B\cap \partial H$, a sphere with four holes.

A \textit{slope disk for $\{\lambda,\rho\}$} is an essential disk in $H$,
possibly separating, which is contained in $B-F$ and is not isotopic to any
component of~$F$. The boundary of a slope disk always separates $\Sigma$
into two pairs of pants, conversely any loop in $\Sigma$ that is not
homotopic into $\partial \Sigma$ is the boundary of a unique slope disk.
(Throughout our work, ``unique'' means unique up to isotopy in an
appropriate sense.) If two slope disks are isotopic in $H$, then they are
isotopic in~$B$.

An arc in $\Sigma$ whose endpoints lie in two different boundary circles of
$\Sigma$ is called a \textit{cabling arc.}  Figure~\ref{fig:slopedisks} shows a
pair of cabling arcs disjoint from a slope disk. A slope disk is disjoint
from a unique pair of cabling arcs, and each cabling arc determines a
unique slope disk.

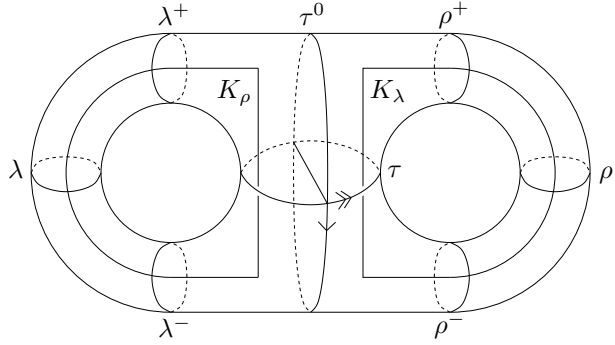
\begin{figure}
\labellist
\small \hair 2pt
\pinlabel $\lambda^+$ [B] at 148 298
\pinlabel $\rho^+$ [B] at 437 300
\pinlabel $\lambda^-$ [B] at 149 -25
\pinlabel $\rho^-$ [B] at 433 -21
\pinlabel $\lambda$ [B] at -15 140
\pinlabel $\tau$ [B] at 376 141
\pinlabel $\rho$ [B] at 594 140
\pinlabel $K_\rho$ [B] at 210 222
\pinlabel $K_\lambda$ [B] at 369 222
\pinlabel $\tau^0$ [B] at 291 298
\endlabellist
\begin{center}
\includegraphics[width=45 ex]{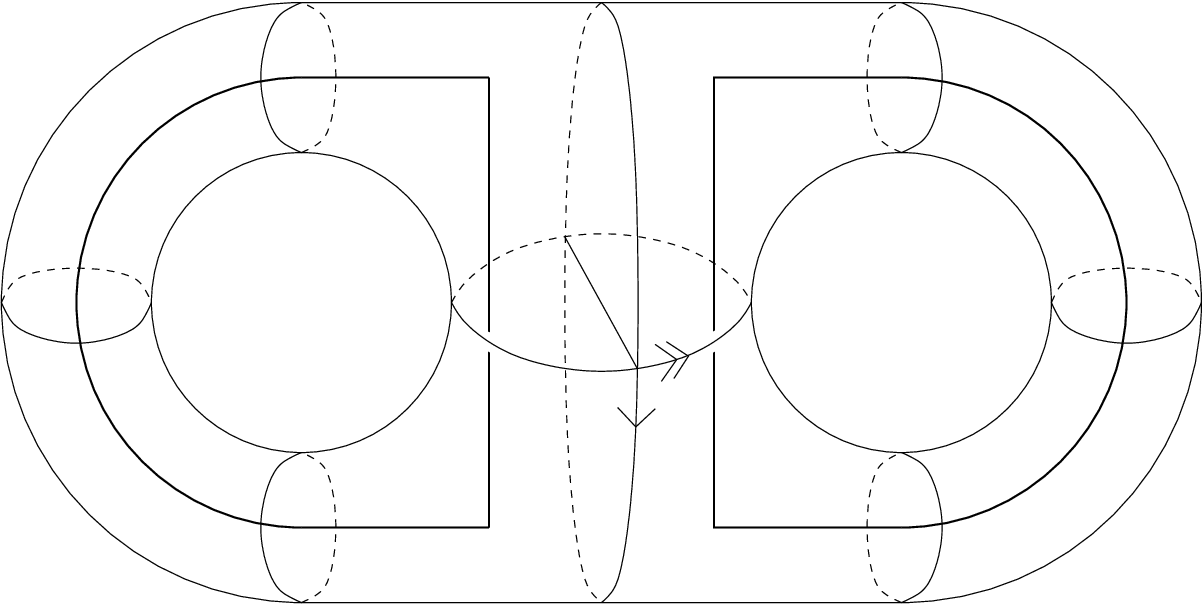}
\caption{The slope-zero perpendicular disk $\tau^0$. It is chosen so that
$K_\lambda$ and $K_\rho$ have linking number~$0$.}
\label{fig:slope_coords}
\end{center}
\end{figure}

Each choice of nonseparating slope disk for a pair $\mu=\{\lambda,\rho\}$
determines a correspondence between~$\Q\cup\{\infty\}$ and the set of all
slope disks of $\mu$, as follows. Fixing a nonseparating slope disk $\tau$
for $\mu$, write $(\mu;\tau)$ for the ordered pair consisting of $\mu$
and~$\tau$.
\begin{definition} A \textit{perpendicular disk} for $(\mu;\tau)$
is a disk $\tau^\perp$, with the following properties:
\begin{enumerate}
\item $\tau^\perp$ is a slope disk for $\mu$.
\item $\tau$ and $\tau^\perp$ intersect transversely in one arc.
\item $\tau^\perp$ separates $H$.
\end{enumerate}
\end{definition}
There are infinitely many choices for $\tau^\perp$, but because $H\subset
S^3$ there is a natural way to choose a particular one, which we call
$\tau^0$. It is illustrated in Figure~\ref{fig:slope_coords}. To construct
it, start with any perpendicular disk and change it by Dehn twists of $H$
about $\tau$ until the core circles of the complementary solid tori have
linking number~$0$ in~$S^3$.

For calculations, it is convenient to draw the picture as in
Figure~\ref{fig:slope_coords}, and orient the boundaries of $\tau$ and
$\tau^0$ so that the orientation of $\tau^0$ (the ``$x$-axis''), followed
by the orientation of $\tau$ (the ``$y$-axis''), followed by the outward
normal of $H$, is a right-hand orientation of $S^3$. At the other
intersection point, these give the left-hand orientation, but the
coordinates are unaffected by changing the choices of which of
$\{\lambda,\rho\}$ is $\lambda$ and which is $\rho$, or changing which of
the disks $\lambda^+$, $\lambda^-$, $\rho^+$, and $\rho^-$ are ``$+$'' and
which are ``$-$'', provided that the ``$+$'' disks both lie on the same
side of $\lambda\cup\rho\cup\tau$ in Figure~\ref{fig:slope_coords}.

\begin{figure}
\labellist
\small \hair 2pt
\pinlabel $\lambda^+$ [B] at 65 66
\pinlabel $\lambda^-$ [B] at 138 66
\pinlabel $\lambda^+$ [B] at 209 66
\pinlabel $\lambda^-$ [B] at 281 66
\pinlabel $\rho^+$ [B] at 65 142
\pinlabel $\rho^-$ [B] at 138 142
\pinlabel $\rho^+$ [B] at 209 142
\pinlabel $\rho^-$ [B] at 281 142
\pinlabel $\lambda^+$ [B] at 65 210
\pinlabel $\lambda^-$ [B] at 138 210
\pinlabel $\lambda^+$ [B] at 209 210
\pinlabel $\lambda^-$ [B] at 281 210
\pinlabel $\rho^+$ [B] at 65 286
\pinlabel $\rho^-$ [B] at 138 286
\pinlabel $\rho^+$ [B] at 209 286
\pinlabel $\rho^-$ [B] at 281 286
\pinlabel $\lambda^+$ [B] at 420 0
\pinlabel $\lambda^-$ [B] at 740 0
\pinlabel $\rho^+$ [B] at 420 310
\pinlabel $\rho^-$ [B] at 745 309
\pinlabel $\tau^0$ [B] at 738 158
\pinlabel $\tau$ [B] at 578 244
\endlabellist
\begin{center}
\includegraphics[width=\textwidth]{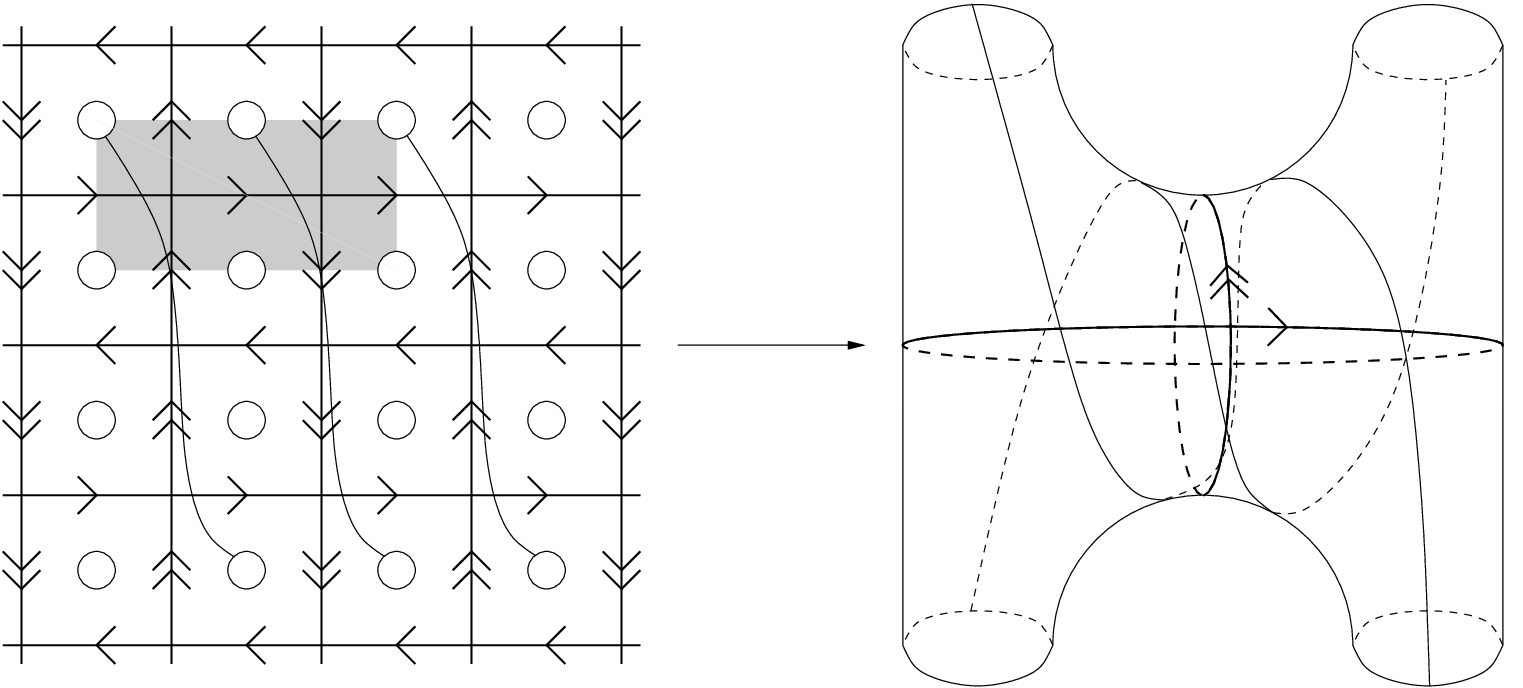}
\caption{The covering space $\widetilde{\Sigma}\to\Sigma$, and some lifts
of a pair of $[1,-3]$-cabling arcs. The shaded region is a fundamental
domain.}
\label{fig:covering}
\end{center}
\end{figure}
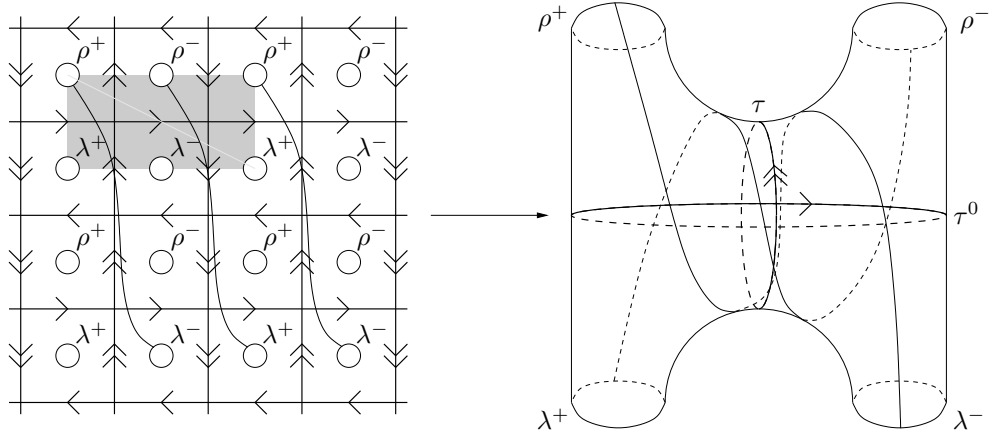

Let $\widetilde{\Sigma}$ be the covering space of $\Sigma$ such that:
\begin{enumerate}
\item $\widetilde{\Sigma}$ is the plane with an open disk of radius $1/8$
removed from each point with coordinates in $\Z\times \Z+(\frac12,\frac12)$.
\item The components of the preimage of $\tau$ are the vertical lines
with integer $x$-coordinate.
\item The components of the preimage of $\tau^0$ are the horizontal lines
with integer $y$-coordinate.
\end{enumerate}
\noindent Figure~\ref{fig:covering} shows a picture of $\widetilde{\Sigma}$
and a fundamental domain for the action of its group of covering
transformations, which is the orientation-preserving subgroup of the group
generated by reflections in the half-integer lattice lines (that pass
through the centers of the missing disks). Each circle of
$\partial\widetilde{\Sigma}$ double covers a circle of~$\partial \Sigma$.

Each lift of a cabling arc $\alpha$ of $\Sigma$ to $\widetilde{\Sigma}$
runs from a boundary circle of $\widetilde{\Sigma}$ to one of its
translates by a vector $(p,q)$ of signed integers, defined up to
multiplication by the scalar $-1$. In this way $\alpha$ receives a
\textit{slope pair} $[p,q]=\{(p,q),(-p,-q)\}$, and is called a
\textit{$[p,q]$-cabling arc.}  The corresponding slope disk is assigned the
slope pair $[p,q]$ as well.

An important observation is that a $[p,q]$-slope disk is nonseparating in
$H$ if and only if $q$ is odd. Both happen exactly when a corresponding
cabling arc has one endpoint in $\lambda^+$ or $\lambda^-$ and the other in
$\rho^+$ or~$\rho^-$.

\begin{definition} Let $\lambda$, $\rho$, and $\tau$ be as above, and let
$\mu=\{\lambda,\rho\}$. The \textit{$(\mu;\tau)$-slope} of a $[p,q]$-slope 
disk or cabling arc is~$q/p\in \Q\cup\{\infty\}$.
\end{definition}
\noindent The $(\mu;\tau)$-slope of $\tau^0$ is $0$, and the
$(\mu;\tau)$-slope of $\tau$ is~$\infty$.

Slope disks for a primitive pair are called \textit{simple} disks, and are
handled in a special way. Rather than using a particular choice of $\tau$
from the context, one chooses $\tau$ to be some third primitive disk.
Altering this choice can change $[p,q]$ to any $[p+nq,q]$, but the quotient
$p/q$ is well-defined as an element of $\Q/\Z\cup\{\infty\}$. This element
$[p/q]$ is called the \textit{simple slope} of the slope disk. The simple
slope is $[0]$ exactly when the slope disk is itself primitive, and has $q$
odd exactly when the simple disk is nonseparating. Simple disks have the
same simple slope exactly when they are equivalent by an element of the
Goeritz group.

\section{The cabling construction}
\label{sec:cabling}

In a sentence, the cabling construction is to ``Think of the union of $K$
and the tunnel arc as a $\theta$-curve, and rationally tangle the ends of
the tunnel arc and one of the arcs of $K$ in a neighborhood of the other
arc of $K$.''  We sometimes call this ``swap and tangle,'' since one of the
arcs in the knot is exchanged for the tunnel arc, then the ends of other
arc of the knot and the tunnel arc are connected by a rational tangle.

\begin{figure}
\labellist
\small \hair 2pt
\pinlabel $\pi_0$ [B] at -7 178
\pinlabel $\pi$ [B] at 65 178
\pinlabel $\pi_1$ [B] at 120 178
\pinlabel $\pi_0$ [B] at 177 178
\pinlabel $\tau_0$ [B] at 239 227
\pinlabel $\pi_1$ [B] at 304 177
\pinlabel $\pi_1$ [B] at 66 121
\pinlabel $\tau_0$ [B] at 127 93
\pinlabel $\pi_0$ [B] at 93 39
\pinlabel $\tau_1$ [B] at 226 54
\pinlabel $\pi_0$ [B] at 250 32
\pinlabel $\tau_0$ [B] at 290 99
\endlabellist 
\begin{center}
\includegraphics[width=50 ex]{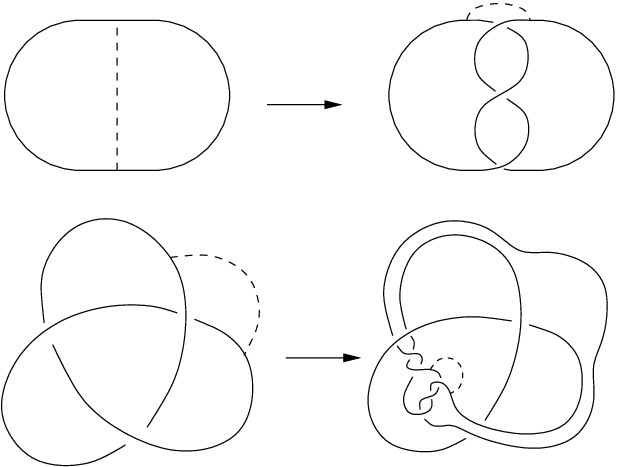}
\caption{Examples of the cabling construction.}
\label{fig:cabling}
\end{center}
\end{figure}
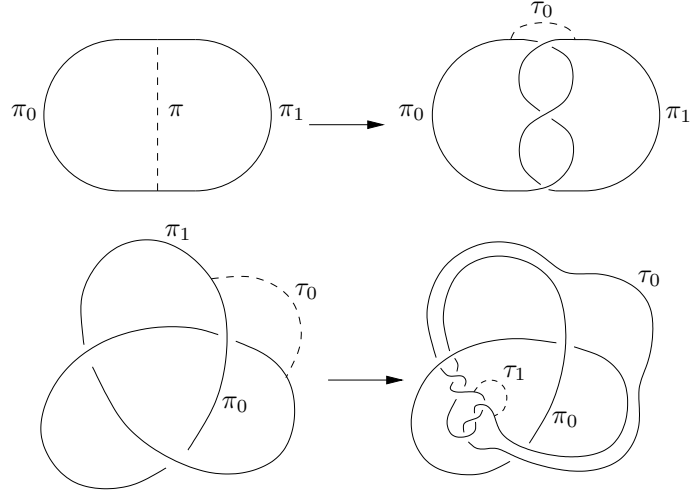

Figure~\ref{fig:cabling} illustrates two cabling constructions, one
starting with the trivial knot and obtaining the trefoil, then another
starting with the tunnel of the trefoil.

More precisely, begin with a triple $\{\lambda,\rho,\tau\}$, regarded as a
pair $\mu=\{\lambda,\rho\}$ with a slope disk $\tau$. Choose one of the
disks in $\{\lambda,\rho\}$, say $\lambda$, and a nonseparating slope disk
$\tau'$ of the pair $\{\lambda,\tau\}$, \textit{other than $\rho$.} This is
a cabling operation producing the tunnel $\tau'$ from $\tau$. In terms of
the ``swap and tangle'' description of a cabling, $\lambda$ is dual to the
arc of $K_\tau$ that is retained, and the slope disk $\tau'$ determines a
pair of cabling arcs that form the rational tangle that replaces the arc of
$K_\tau$ dual to~$\rho$.

Provided that $\{\lambda,\rho,\tau\}$ was not a primitive triple, we define
the \textit{slope} of this cabling operation to be the
$(\{\lambda,\tau\};\rho)$-slope of~$\tau'$.  When $\{\lambda,\rho,\tau\}$
is primitive, the cabling construction starts with the tunnel of the
trivial knot and produces an upper or lower tunnel of a $2$-bridge knot,
unless $\tau'$ is primitive, in which case it is again the tunnel of the
trivial knot and the cabling is called \textit{trivial.} The slope of a
cabling starting with a primitive triple is defined to be the simple slope
of $\tau'$. The cabling is trivial when the simple slope is~$[0]$.

Theorem~13.2 of~\cite{CM} shows that every tunnel of every tunnel
number~$1$ knot can be obtained by a uniquely determined sequence of cabling
constructions. The associated cabling slopes form a sequence
\[ m_0,\;m_1,\;\cdots\;,\;m_n = [p_0/q_0],\;q_1/p_1,\;\cdots\;,\;q_n/p_n\]
where $m_0\in\Q/\Z$ and each $q_i$ is odd. 

There is a second set of invariants associated to a tunnel. Each $m_i$ is
the slope of a cabling that begins with a triple of disks
$\{\lambda_{i-1},\rho_{i-1},\tau_{i-1}\}$ and finishes with
$\{\lambda_i,\rho_i,\tau_i\}$. For $i\geq 2$, put $s_i=1$ if
$\{\lambda_i,\rho_i,\tau_i\}=\{\tau_{i-2},\tau_{i-1},\tau_i\}$, and $s_i=0$
otherwise. In terms of the swap-and-tangle construction, the invariant
$s_i$ is $1$ exactly when the rational tangle replaces the arc that was
retained by the previous cabling (for $i=1$, the choice does not matter, as
there is an element of the Goeritz group that preserves $\tau_0$ and
interchanges $\lambda_0$ and $\rho_0$).

In the sequence of triples described in the previous paragraph, the disks
$\lambda_i$ and $\rho_i$ form the \textit{principal pair} for the tunnel
$\tau_i$. They are the disks called $\mu^+$ and $\mu^-$
in~\cite{Scharlemann-Thompson}.

A nontrivial tunnel $\tau_0$ produced from the tunnel of the trivial knot
by a single cabling construction is called a \textit{simple} tunnel. As
already noted, these are the ``upper and lower'' tunnels of $2$-bridge
knots. Not surprisingly, the simple slope $m_0$ is a version of the
standard rational parameter that classifies the $2$-bridge
knot~$K_{\tau_0}$.

A tunnel is called \textit{semisimple} if it is disjoint from a primitive
disk, but not from any primitive pair. The simple and semisimple tunnels
are exactly the $(1,1)$-tunnels, that is, the upper and lower tunnels of
knots in $1$-bridge position with respect to a Heegaard torus of~$S^3$.
A tunnel is semisimple if and only if all $s_i=0$. The reason is that both
conditions characterize cabling sequences in which one of the original
primitive disks is retained in every cabling; this corresponds to the fact
that the union of the tunnel arc and one of the arcs of the knot is
unknotted. 

A tunnel is called \textit{regular} if it is neither primitive, simple, or
semisimple.

\section{The middle tunnels}
\label{sec:regular}

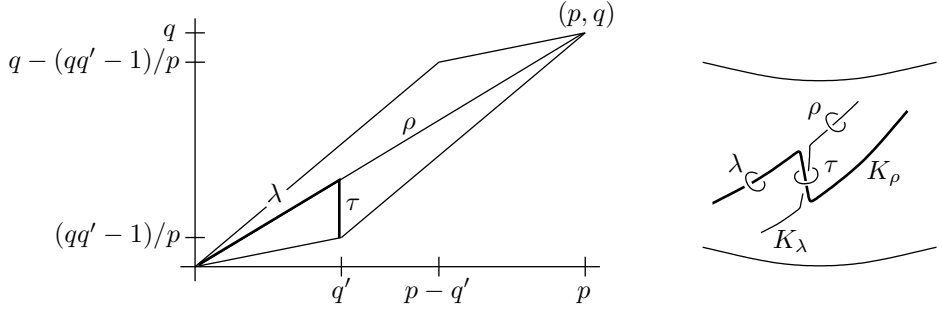
\begin{figure}
\labellist
\small \hair 2pt
\pinlabel $q$ [B] at -3 85
\pinlabel $q-(qq'-1)/p$ [B] at -29 74
\pinlabel $(qq'-1)/p$ [B] at -21 14
\pinlabel $\lambda$ [B] at 33 28
\pinlabel $\tau$ [B] at 60 25
\pinlabel $\rho$ [B] at 79 53
\pinlabel $q'$ [B] at 56 -5
\pinlabel $p-q'$ [B] at 89 -5
\pinlabel $p$ [B] at 139 -5
\pinlabel $(p,q)$ [B] at 140 90
\pinlabel $\lambda$ [B] at 191 38
\pinlabel $\rho$ [B] at 218 58
\pinlabel $\tau$ [B] at 224 38
\pinlabel $K_{\lambda}$ [B] at 210 12
\pinlabel $K_{\rho}$ [B] at 242 36
\endlabellist 
\begin{center}
\includegraphics[width=0.8\textwidth]{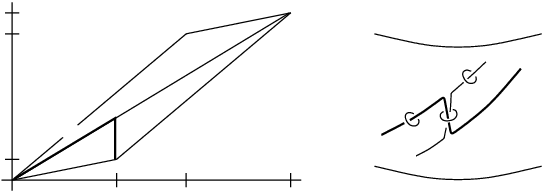}
\caption{The properties of $q'$. The darker segments correspond to
$K_\rho$, a $(q',(qq'-1)/p)$ torus knot. The picture on the right shows
$K_\lambda$ in the torus $T\subset S^3$, and $K_\rho$ pulled slightly
outside of $T$.}
\label{fig:qprime}
\end{center}
\end{figure}
In this section we will prove Theorem~\ref{thm:middle_tunnels_slopes}. We
have relatively prime integers $p>q\geq 2$, and we use the notations $T$,
$V$, and $W$ of Section~\ref{sec:statement}.

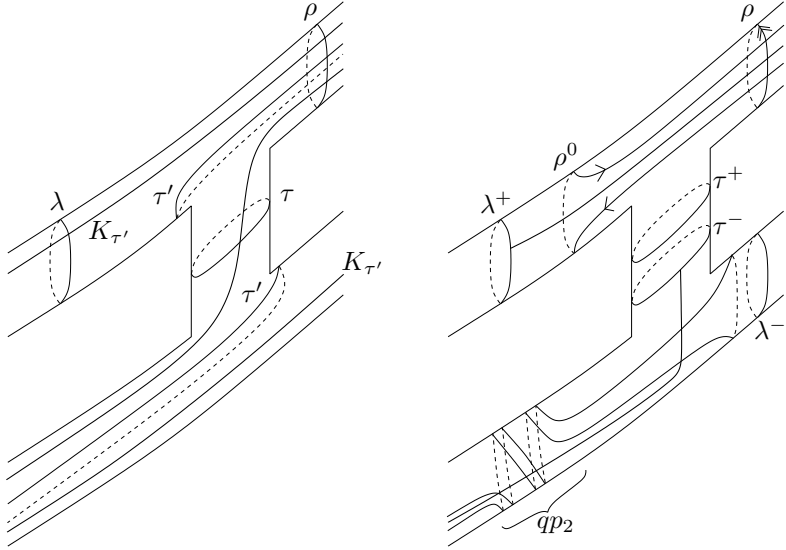
\begin{figure}
\labellist
\small \hair 2pt
\pinlabel $\lambda$ [B] at 45 290
\pinlabel $\rho$ [B] at 263 458
\pinlabel $K_{\tau'}$ [B] at 90 265
\pinlabel $\tau'$ [B] at 137 295
\pinlabel $\tau$ [B] at 243 298
\pinlabel $\tau'$ [B] at 213 211
\pinlabel $K_{\tau'}$ [B] at 308 238
\pinlabel $\lambda^+$ [B] at 423 288
\pinlabel $\rho^0$ [B] at 483 329
\pinlabel $\rho$ [B] at 639 457
\pinlabel $\tau^+$ [B] at 623 309
\pinlabel $\tau^-$ [B] at 623 270
\pinlabel $qp_2$ [B] at 474 16
\pinlabel $\lambda^-$ [B] at 659 179
\endlabellist 
\begin{center}
\includegraphics[width=0.82\textwidth]{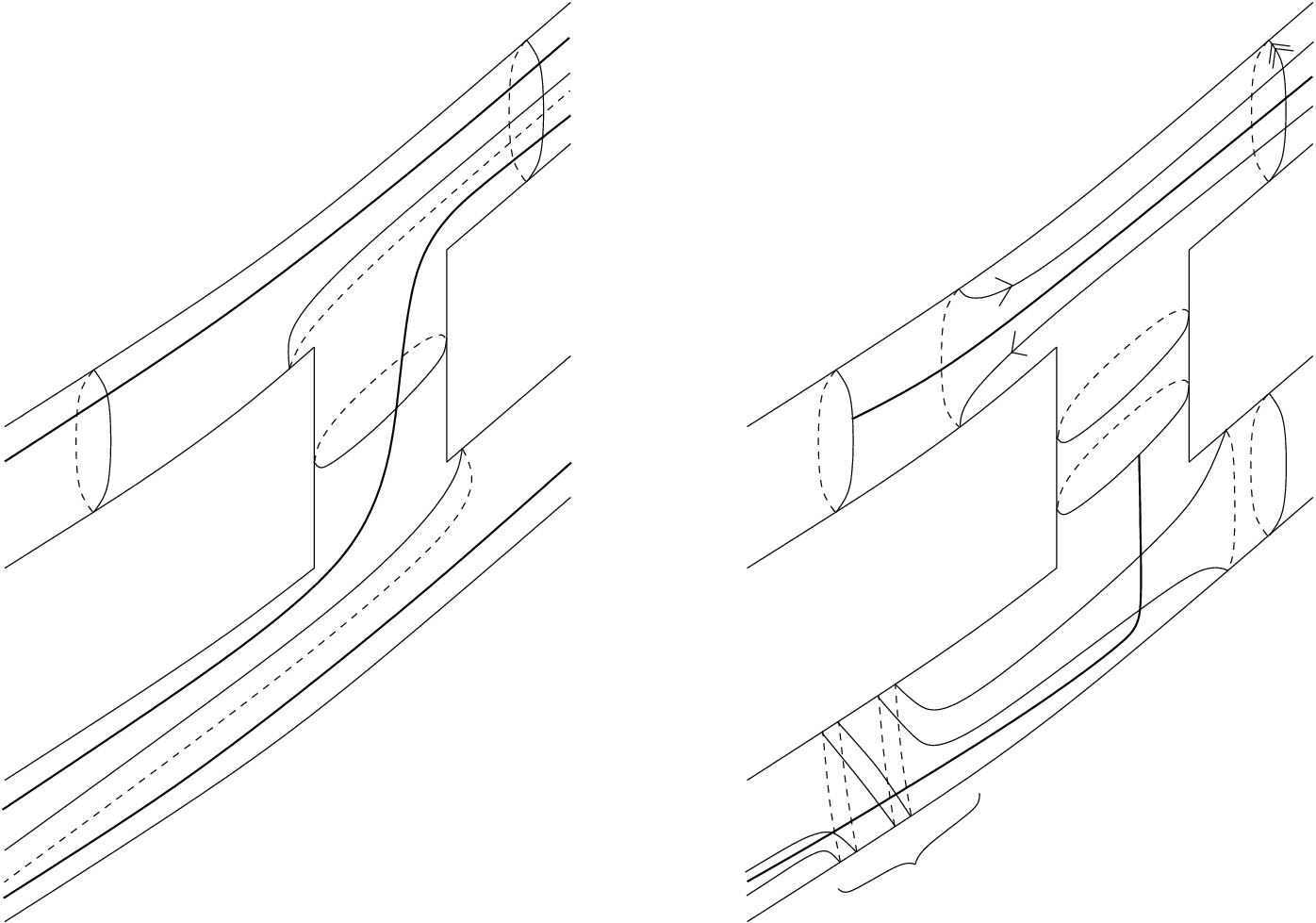}
\caption{The cabling construction that replaces $\rho$ (compare with
Figure~\ref{fig:slope_coords}) with $\tau'$. The left drawing shows the new
tunnel disk $\tau'$ and the knot $K_{\tau'}$. The right drawing shows a
cabling arc for $\tau'$, running from $\lambda^+$ to $\tau^-$, and the
disks $\rho$ and $\rho^0$ used to calculate its slope.}
\label{fig:torus_tunnel}
\end{center}
\end{figure}

First we examine a cabling operation that takes the middle tunnel $\tau$
and produces a middle tunnel of a new torus knot. Let $q'$ be the integer
with $0<q'<p$ such that $qq'\equiv 1\pmod{p}$. If the principal pair
$\{\lambda,\rho\}$ of $\tau$ is positioned as shown in
Figures~\ref{fig:qprime} and~\ref{fig:torus_tunnel} (our inductive
construction of these tunnels will show that the pair shown in the figures
is indeed the principal pair), then $K_\rho$ is a $(q',(qq'-1)/p)$ torus
knot, and $K_\lambda$ is a $(p-q', q-(qq'-1)/p)$ torus knot. We set
$(p_1,q_1)= (q',(qq'-1)/p)$ and $(p_2,q_2) = (p-q', q-(qq'-1)/p)$, so that
$K_\rho$ and $K_\lambda$ are respectively the $(p_1,q_1)$ and $(p_2,q_2)$
torus knots.

In Figure~\ref{fig:qprime}, the linking number of $K_\rho$ with
$K_\lambda$, up to sign conventions, is $q_1p_2$. One way to see this is to
note that a Seifert surface for $K_\lambda$ can be constructed using $q_2$
meridian disks of $V$ and $p_2$ meridian disks of $W$ (by attaching bands
contained in a small neighborhood of $T$). When $K_\rho$ is pulled slightly
outside of $V$, as indicated in Figure~\ref{fig:qprime}, it has $q_1$
intersections with each of the $p_2$ meridian disks of $W$, all crossing
the disks in the same direction.

Figure~\ref{fig:torus_tunnel} shows the new tunnel disk $\tau'$ for a
cabling construction that produces a $(p+p_2,q+q_2)$ torus knot
$K_{\tau'}$. This disk meets $T$ perpendicularly. The drawing on the right
in Figure~\ref{fig:torus_tunnel} illustrates the setup for the calculation
of the $(\{\lambda,\tau\};\rho)$-slope pair of $\tau'$. The $qp_2$ turns of
$\rho^0$, with the case $qp_2=2$ drawn in the figure, make the copies of
$K_\tau=K_{p,q}$ and $K_\lambda=K_{p_2,q_2}$ in its complement have linking
number $0$. A cabling arc for $\tau'$ is shown. Examination of its
crossings with $\partial \rho$ and $\partial \rho^0$ shows that the slope
pair of $\tau'$ is~$[1,2qp_2+1]$.

Put $U=\begin{pmatrix}1&1\\0&1\end{pmatrix}$ and
$L=\begin{pmatrix}1&0\\1&1\end{pmatrix}$. If $K_1$ is a $(p_1,q_1)$ torus
knot and $K_2$ is a $(p_2,q_2)$ torus knot, we denote by $M(K_1,K_2)$ the
matrix $\begin{pmatrix} p_1 & q_1 \\ p_2 & q_2\end{pmatrix}$. In our case,
this is the matrix $M(K_\rho, K_\lambda)$. Adding the rows of
$M(K_\rho,K_\lambda)$ gives $(p,q)$, corresponding to $K_\tau$, so
\[ M(K_\tau,K_\lambda) = U\cdot M(K_\rho,K_\lambda)\ .\]
The left drawing of Figure~\ref{fig:torus_tunnel} can be repositioned by
isotopy so that $\lambda$, $\tau$, and $\tau'$ look respectively as did
$\lambda$, $\rho$, and $\tau$ in the original picture, with $\tau'$ as the
tunnel of the $(p+p_2,q+q_2)$ torus knot. Thus the procedure can be
repeated, each time multiplying the matrix by another factor of~$U$.

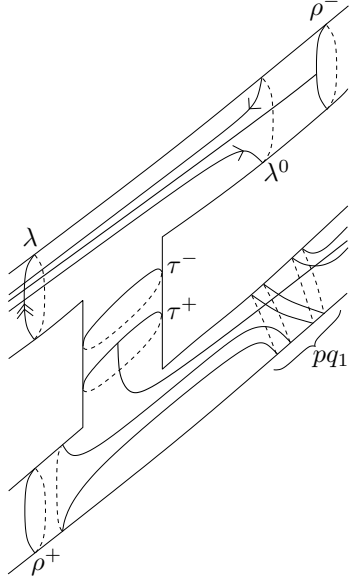
\begin{figure}
\labellist
\small \hair 2pt
\pinlabel $\lambda$ [B] at 20 280
\pinlabel $\rho^-$ [B] at 271 476
\pinlabel $\lambda^0$ [B] at 230 335
\pinlabel $\tau^-$ [B] at 150 257
\pinlabel $\tau^+$ [B] at 150 218
\pinlabel $\rho^+$ [B] at 34 4
\pinlabel $pq_1$ [B] at 276 181
\endlabellist 
\begin{center}
\includegraphics[width=0.36 \textwidth]{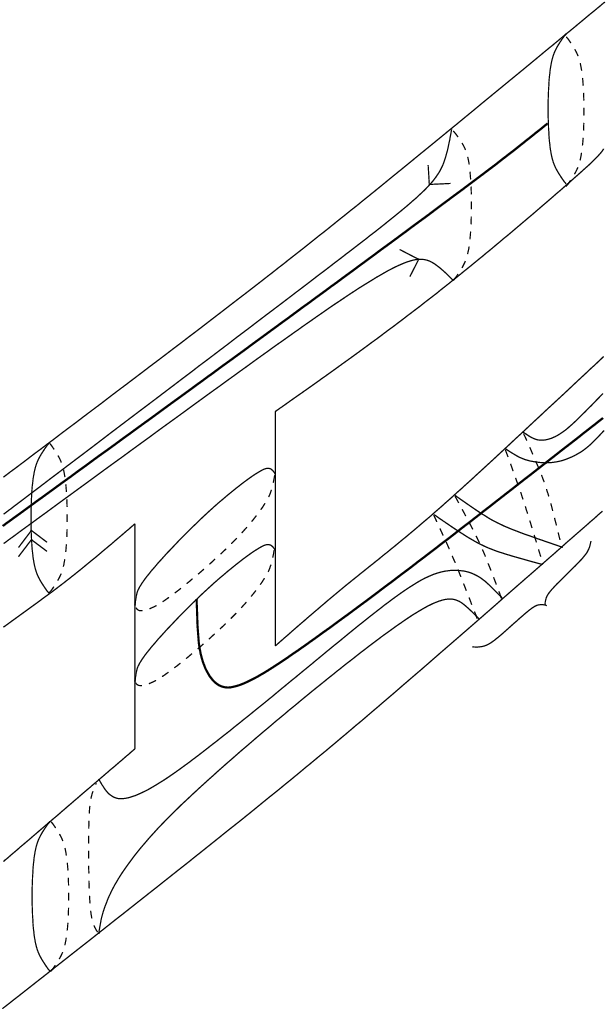}
\caption{Calculation of the slope of $\tau'$ for a cabling construction
replacing $\lambda$. The cabling arc runs from $\rho^-$ to $\tau^+$.}
\label{fig:torus_tunnel1}
\end{center}
\end{figure}
Figure~\ref{fig:torus_tunnel1} illustrates the similar
calculation of the slope of the
cabling construction replacing $\lambda$ by a new tunnel $\tau'$.
This produces a $(p+p_1,q+q_1)$ torus knot. In this case we have
\[ M(K_\rho, K_\tau) = L\cdot M(K_\rho,K_\lambda)\ .\]
The slope pair of $\tau'$ is $[1,2qp_1-1]$. One might expect $2qp_1+1$ as
the second term, in analogy with the construction replacing
$\rho$. However, as seen in Figure~\ref{fig:torus_tunnel1}, the $pq_1$
twists needed in $\lambda^0$ are in the same direction as the twists in the
calculation for $\rho$, not in the mirror-image sense. This results in two
fewer crossings of the cabling arc for $\tau'$ with $\lambda^0$ than
before. In fact, the slope pairs for the two constructions can be described
in a uniform way: For either of the matrices $M(K_\tau,K_\lambda)$ and
$M(K_\rho,K_\tau)$, a little bit of arithmetic shows that the second entry
of the slope pair for the cabling operation that produced them is the sum
of the product of the diagonal entries and the product of the off-diagonal
entries, that is, $[1,pq_2+qp_2]$ in the first case and $[1,pq_1+qp_1]$ in
the second.

We can now describe the complete cabling sequence. Still assuming that $p$
and $q$ are both positive and $p>q$, write $p/q$ as $[n_1,n_2,\ldots, n_k]$
with all $n_i$ positive. We may assume that $n_k\neq 1$. According as $k$
is even or odd, we consider the product $U^{n_k}L^{n_{k-1}}\cdots
U^{n_2}L^{n_1}$ or $L^{n_k}U^{n_{k-1}}\cdots U^{n_2}L^{n_1}$. 

Start with a trivial knot regarded as a $(1,1)$ torus knot, and the
``middle'' tunnel $\tau$ in $T$. For the disks $\lambda$ and $\rho$ shown
in Figure~\ref{fig:torus_tunnel}, $K_\rho$ is a $(1,0)$ torus knot and
$K_\lambda$ is a $(0,1)$ torus knot. For this positioning of the trivial
knot, the disks $\lambda$, $\rho$, and $\tau$ are all primitive, so
$\{\lambda,\rho\}$ may be regarded as the principal pair for the tunnel
$\tau$. Cablings of the two types above will preserve the fact that
the pair $\{\lambda,\rho\}$ shown in Figure~\ref{fig:torus_tunnel} is the
principal pair.

At this point, the matrix $M(K_\rho,K_\lambda)$ is the identity matrix.
Multiplying by $L^{n_1}$ corresponds to doing $n_1$ cabling constructions
of the second type described above (replacing $\lambda$). These cablings
have slope $[1]=[0]\in \Q/\Z$, so are trivial cablings, but their effect is
to produce the trivial knot positioned as an $(n_1+1,1)$ torus knot. At
that stage, $M(K_\rho,K_\lambda)$ is $\begin{pmatrix} 1 & 0 \\ n_1 &
1\end{pmatrix}$, which is the matrix $M_{-1}$ in the statement of
Theorem~\ref{thm:middle_tunnels_slopes}. Then, multiplying by $U$
corresponds to a nontrivial cabling construction of the first type
(replacing $\rho$). The new matrix $M(K_\tau,K_\lambda)$ is
$\begin{pmatrix} n_1+1 & 1 \\ n_1 & 1\end{pmatrix}$, or in the notation of
Theorem~\ref{thm:middle_tunnels_slopes}, $M_0$, and the knot $K_{\tau'}$ is
a $(2n_1+1,2)$ torus knot. As explained above, the construction has slope
pair $[1,2n_1+1]$, so the simple slope is $m_0=[1/(2n_1+1)]$. Continue by
multiplying $n_2-1$ additional times by $U$, then $n_3$ times by $L$ and so
on, performing additional cabling constructions with slopes calculated as
above from the matrices of the current $K_\rho$, $K_\lambda$, and $K_\tau$.
This produces the sequence of matrices $M_t$ in
Theorem~\ref{thm:middle_tunnels_slopes} and the corresponding slope
invariants $m_t=a_td_t+b_tc_t$.

At the end, there is no cabling construction corresponding to the last
factor $L$ or $U$. For specificity, suppose $k$ was even and the product
was $U^{n_k}L^{n_{k-1}}\cdots U^{n_1}L^{n_1}$.  At the last stage, we apply
$n_k-1$ cabling constructions corresponding to multiplications by $U$, and
arrive at a tunnel $\tau$ for which $M(K_\rho,K_\lambda)$ is
$U^{n_k-1}L^{n_{k-1}}\cdots U^{n_2}L^{n_1}=M_N$. The sum of the rows is
then $(p,q)$ (multiplying by $U$ and using the case ``$q/s$'' of Lemma~14.3
of~\cite{CM}), so $K_\tau$ is the $(p,q)$ torus knot. The case when $k$ is
odd is similar (multiplying by $L$ and using the ``$p/r$'' case of
Lemma~14.3 of~\cite{CM}). In summary, there are $-1+\sum_{i=2}^{k} n_i=N+1$
nontrivial cabling constructions, whose slopes can be calculated as
in Theorem~\ref{thm:middle_tunnels_slopes}.

When $m_t$ is calculated from the matrix
$\begin{pmatrix} a_t & b_t\\ c_t & d_t\end{pmatrix}$,
the knot $K_\tau$ is an $(a_t+c_t,b_t+d_t)$ torus knot, which is part~(ii)
of Theorem~\ref{thm:middle_tunnels_slopes}. For part~(iii), we have $s_t=1$
when the constructions change from replacing $\rho$ to replacing $\lambda$,
or vice versa. This occurs when we change from multiplying by $U$ to
multiplying by $L$, or vice versa, that is, when $A_t\neq A_{t-1}$.

\section{The upper and lower tunnels}
\label{sec:semisimple}

\begin{figure}
\begin{center}
\includegraphics[width=0.8\textwidth]{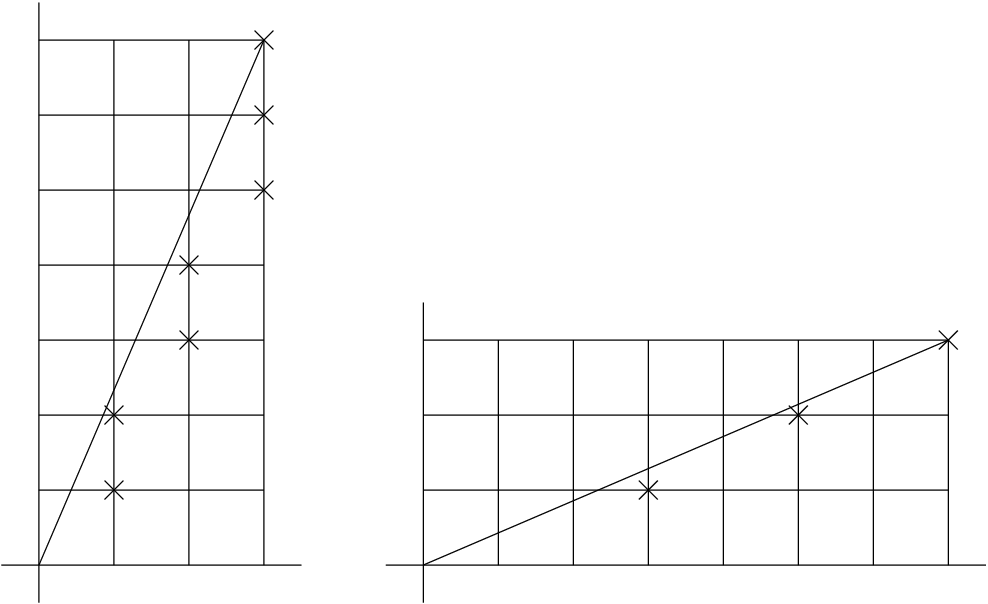}
\caption{The points $(p_k,k)$ for $1\leq k\leq q$ for
the cases $(p,q)=(3,7)$ and $(p,q)=(7,3)$.}
\label{fig:pk}
\end{center}
\end{figure}

Again we use the notations $T$, $V$, and $W$ of previous sections. Denoting
the unit interval $[0,1]$ by $I$, fix a product $T\times I \subset W$ with
$T=T\times\{0\}$. 
\begin{definition} Let $p$ and $q$ be relatively prime integers, both
greater than $1$. For integers $k$ with
$1\leq k\leq q$, put
\[p_k = \lceil kp/q\rceil = \min\{ j\;|\; jq/p \geq
k\}\ .\]
\end{definition}
\noindent Figure~\ref{fig:pk} shows the points $(p_k,k)$ for $1\leq k\leq q$ 
for the cases $(p,q)=(3,7)$ and $(p,q)=(7,3)$.

\begin{figure}
\begin{center}
\includegraphics[width=0.8\textwidth]{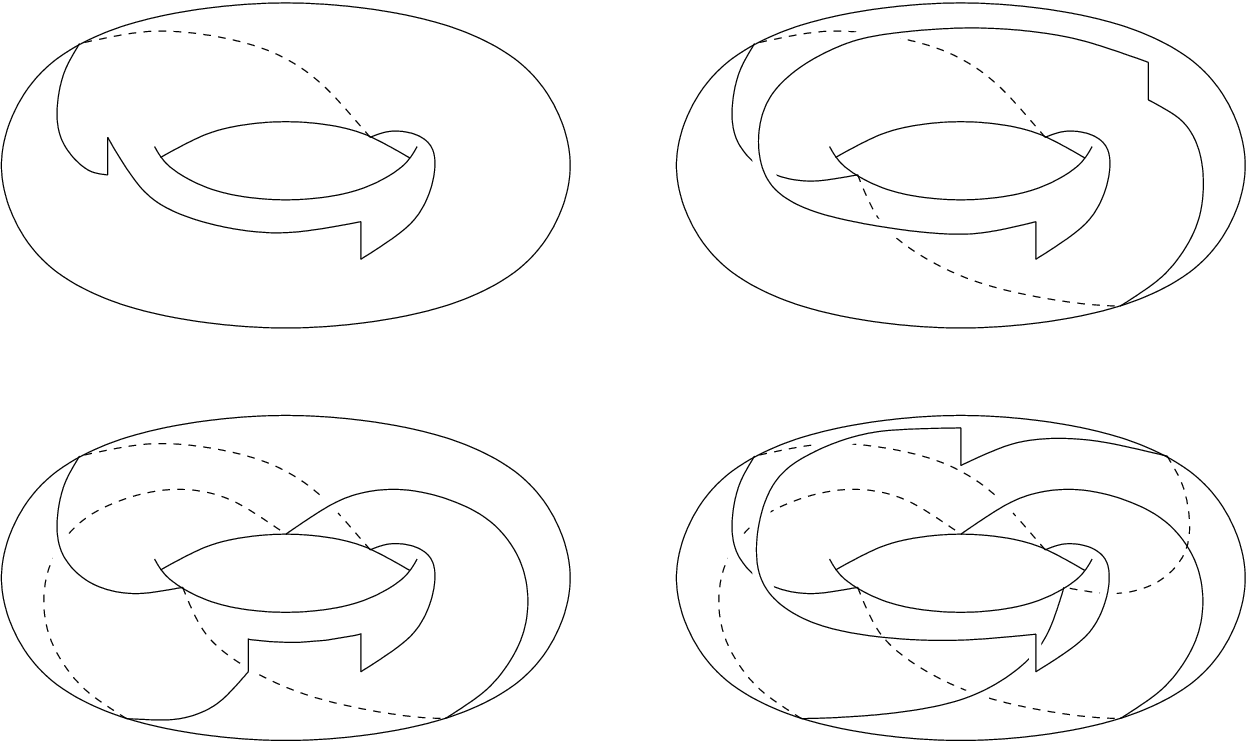}
\caption{The knots $K(3,5;1)$, $K(3,5;2)$, $K(3,5;3)$, and $K(3,5;4)$.
The first two are trivial, $K(3,5;3)$ is a $(2,3)$ torus knot, and
$K(3,5;4)$ is a $(2,5)$ torus knot.}
\label{fig:ks}
\end{center}
\end{figure}
Define knots $K(p,q;k)\subset T\times I$ as follows.
\begin{definition}
In the universal cover $\R^2\times I$ of $T\times I$, take the arc
(i.~e.~line segment) from $((0,0),0)$ to $((kp/q,k),0)$. If $k<q$, add to
this arc the arc from $((kp/q,k),0)$ to $((kp/q,k),1)$, followed by the arc
from $((kp/q,k),1)$ to $((p_k,k),1)$, followed by the arc from
$((p_k,k),1)$ to $((p_k,k),0)$. The image of these arcs in $T\times
I\subset S^3$ is $K(p,q;k)$. In particular, $K(p,q;q)$ is the standard
$(p,q)$ torus knot.\par
\label{def:Kqpk}
\end{definition}
\noindent Figure~\ref{fig:ks} shows the knots $K(3,5;1)$, $K(3,5;2)$,
$K(3,5;3)$, and $K(3,5;4)$.
\medskip

\begin{figure}
\begin{center}
\includegraphics[width=0.8\textwidth]{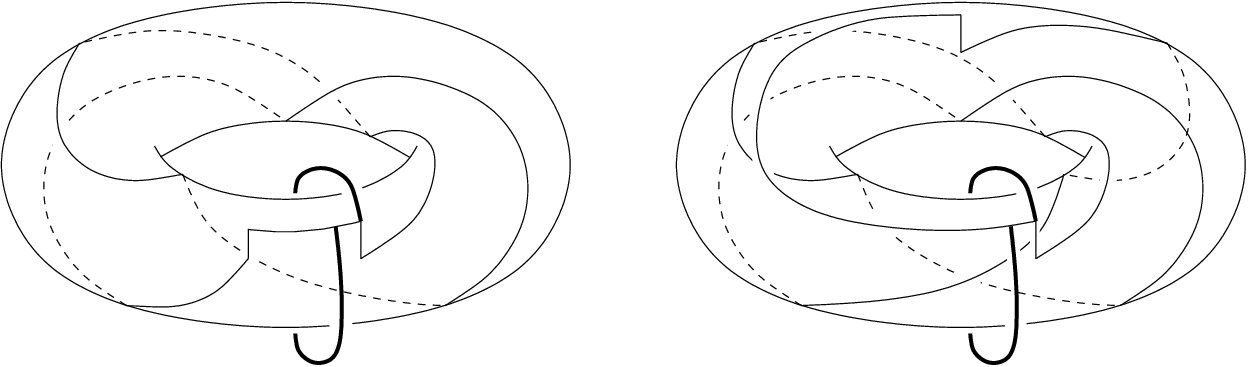}
\caption{Representative arcs of the tunnels $\tau(3,5;3)$ and
$\tau(3,5;4)$.}
\label{fig:upper}
\end{center}
\end{figure}

The \textit{upper tunnel} $\tau(p,q;k)$ of $K(p,q;k)$ is best described by
a picture, given as Figure~\ref{fig:upper}. In particular, $\tau(p, q; q)$
is the standard upper tunnel of the $(p,q)$ torus knot.
Figure~\ref{fig:upper} shows tunnel arcs for the upper tunnels
$\tau(3,5;3)$ and $\tau(3,5;4)$. We will see, inductively, that the unions
of such knots with the particular arcs shown in Figure~\ref{fig:upper} are
the $\theta$-curves dual to the disks of the principal vertex of the tunnel.

\begin{figure}
\begin{center}
\includegraphics[width=0.38\textwidth]{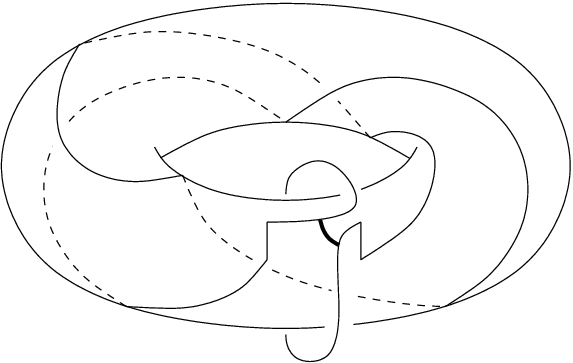}
\caption{The result of the cabling construction producing $\tau(3,5;4)$
from $\tau(3,5;3)$. After isotopy, this becomes the second drawing in
Figure~\ref{fig:upper}.}
\label{fig:uppercabling}
\end{center}
\end{figure}
The cabling construction that takes $\tau(p,q;k)$ to $\tau(p,q;k+1)$ is
illustrated in Figure~\ref{fig:uppercabling} for the case of
$\tau(3,5;3)$. The resulting knot is isotopic to the $K(3,5;4)$ shown in
Figure~\ref{fig:upper}, by pushing the arc that was the tunnel arc of
$\tau(3,5;3)$ down into $T$ and stretching out the new tunnel arc until it
looks like the one in Figure~\ref{fig:upper}.

\begin{figure}
\labellist
\small \hair 2pt
\pinlabel $\lambda^+$ [B] at 198 640
\pinlabel $\tau^+$ [B] at 187 510
\pinlabel $\rho^0$ [B] at 255 548
\pinlabel $\rho$ [B] at 397 735
\pinlabel $\tau^-$ [B] at 622 626
\pinlabel $\lambda^-$ [B] at 619 512
\pinlabel $\tau^+$ [B] at 260 300
\pinlabel $\tau$ [B] at 90 140
\pinlabel $\tau^-$ [B] at 255 -27
\pinlabel $\rho^0$ [B] at 398 300
\pinlabel $\rho$ [B] at 485 141
\pinlabel $\lambda^+$ [B] at 550 300
\pinlabel $\lambda$ [B] at 703 140
\pinlabel $\lambda^-$ [B] at 544 -27
\endlabellist
\begin{center}
\includegraphics[width=0.7\textwidth]{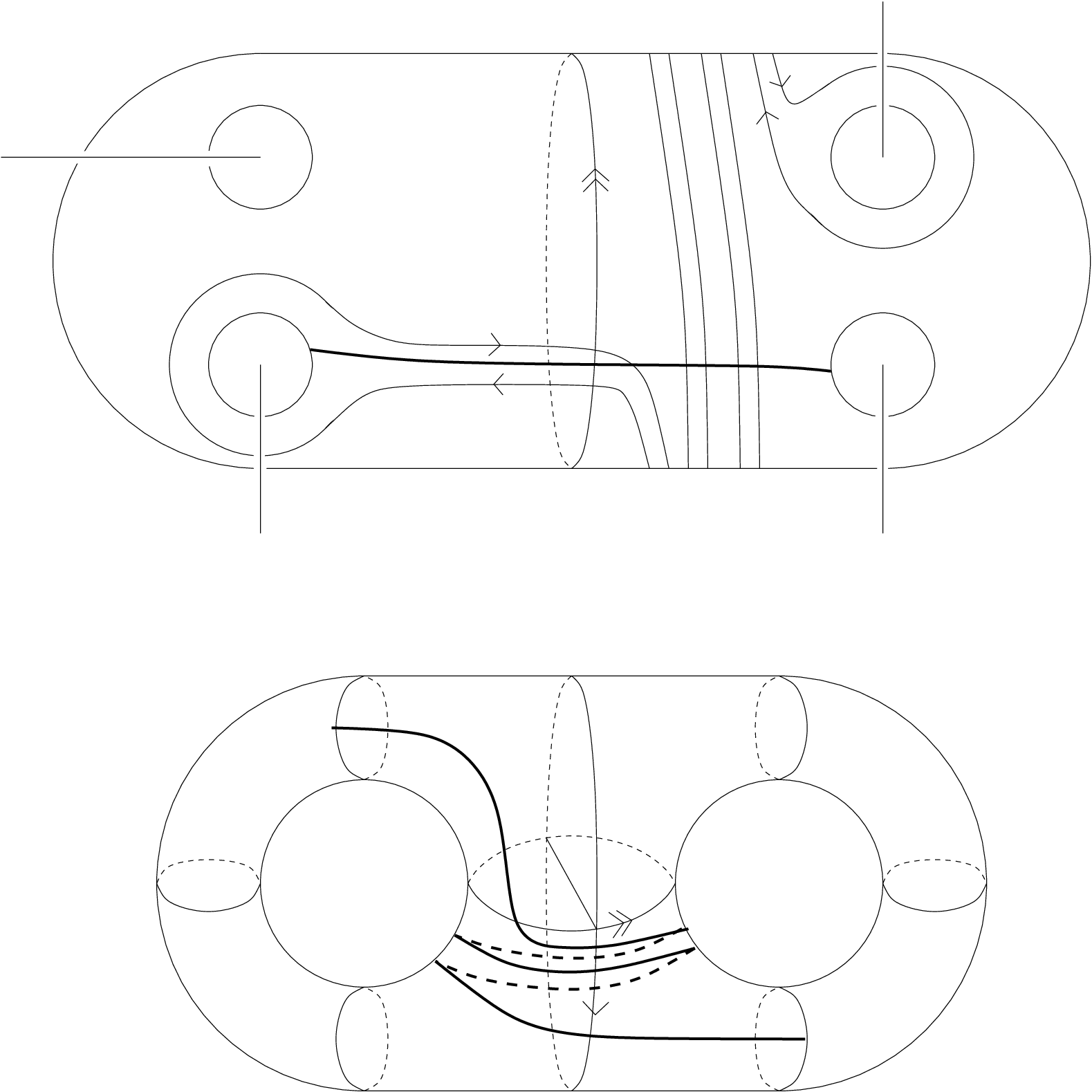}
\caption{The slope calculation for the cabling taking $\tau(p,q;k)$ to
$\tau(p,q,k+1)$. In the top picture, $\rho^0$ makes $p_k$ turns around the
ball; for the example drawn here, $p_k=3$. The cabling arc for the new
tunnel disk has $2p_k-1$ crossings with $\rho^0$, producing the slope pair
$[1,2p_k-1]$ as seen in the standard picture.}
\label{fig:slopes}
\end{center}
\end{figure}
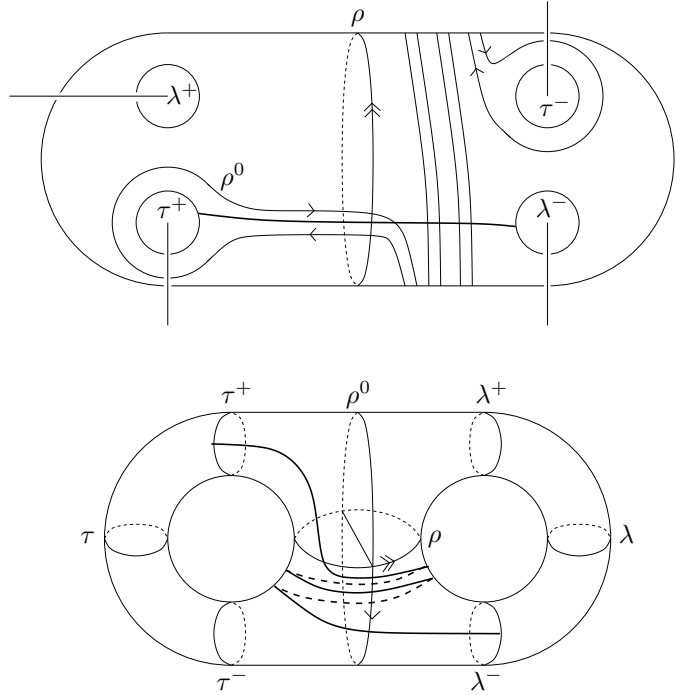
We will now compute the slopes of these cabling operations.
Figure~\ref{fig:slopes} illustrates the calculation for the cabling taking
$\tau=\tau(p,q;k)$ to $\tau(p,q,k+1)$.  The ball shown in the top drawing
in Figure~\ref{fig:slopes} is a regular neighborhood of the arc in the
raised part of $K(p,q;k)$ that connects the endpoints of $\tau(p,q;k)$.
The disk $\rho$ will be replaced.

The $0$-slope disk $\rho^0$ makes $p_k$ turns around the ball.  To see
this, consider a perpendicular disk $\rho^\perp$ for $\rho$ constructed as
follows. In the boundary of the ball in Figure~\ref{fig:slopes}, take an
arc $\alpha$ connecting $\partial \tau^+$ to $\partial \tau^-$, running
across the front of the ball between $\partial \lambda^+$ and $\partial
\lambda^-$, and cutting across $\partial \rho$ in a single point. The
frontier of a regular neighborhood of $\tau^+\cup \alpha\cup\tau^-$ in the
ball is $\rho^\perp$. That is, $\rho^\perp$ is like $\rho^0$ except that it
has no turns around the back of the ball. The representative of $K_\tau$
disjoint from $\rho^\perp$ is isotopic to $K(p,q;k)$, while the
representative of $K_\lambda$ is a core circle of $W$ that completely
encircles this $K(p,q;k)$. In the homology of $V\cup T\times I$, $K(p,q;k)$
represents $p_k$ times the generator, so (for some choice of linking
conventions) $K_\lambda$ has linking number $p_k$ with $K_\tau$. Adding
$p_k$ turns around the ball to $\rho^\perp$ as in the top drawing of
Figure~\ref{fig:slopes} decreases this linking number to $0$, and gives the
perpendicular disk shown in Figure~\ref{fig:slopes}, which must therefore
be~$\rho^0$.

Both diagrams in Figure~\ref{fig:slopes} show the cabling arc for the slope
disk that defines $\tau(p,q;k+1)$, and the bottom picture verifies that its
slope coordinates are $[1,2p_k-1]$. This yields the value for $m_k$ given
in Theorem~\ref{thm:cabling_sequence}. 

We can begin the process with the knot $K(p,q;1)$. For the standard tunnel
arc, all three of the disks $\lambda$, $\rho$, and $\tau$ in the first
drawing of Figure~\ref{fig:slopes} are primitive, since $K_\lambda$,
$K_\rho$, and $K_\tau$ are trivial knots. For $k<k_0$, $p_k=1$ and
$K(p,q;k+1)$ is a trivial knot. This can be seen geometrically, but also
follows inductively from the fact that these cablings have simple slope
$\Big[\displaystyle\frac{1}{2\cdot 1-1}\Big] = [0]\in\Q/\Z$. The process
terminates with the cabling corresponding to $k=p_{q-1}$, which produces
$K(p,q;q)=K_{p,q}$.

\section{Applications}
\label{sec:applications}

Here we will recover the classification of the torus knot tunnels of
M. Boileau, M. Rost, and H. Zieschang~\cite{B-R-Z} and
Y. Moriah~\cite{Moriah}, although not their result that these are all the
tunnels. We consider three cases for $K_{p,q}$:\par
\smallskip

\noindent\textsl{Case I.} $|p-q|=1$.
\smallskip

We may assume that $(p,q)=(n+1,n)$ with $n\geq 2$. For both the upper and
lower tunnels, Theorem~\ref{thm:cabling_sequence} gives $[1/3]$, $5$,
$7,\ldots\,$, $2n-1$ as the slope sequence. For the middle tunnel,
Theorem~\ref{thm:middle_tunnels_slopes} gives the same slope sequence, and
all $s_i=0$, showing that all three tunnels are the same.
\smallskip

\noindent\textsl{Case II.} $|p-q|\neq1$, but $p\equiv \pm1\mod{q}$ or
$q\equiv \pm1\mod{p}$
\smallskip

Again we assume that $p,q\geq 2$, and reduce to the case when $p>q$.
Suppose first that $p=mq+1$ with $m\geq 2$. For the upper tunnel,
Theorem~\ref{thm:cabling_sequence} gives slopes $[1/(2m+1)]$, $4m+1$,
$6m+1,\ldots\,$, $2m(q-1)+1$ (to find the $p_k$, notice that the line
segment in $\R^2$ from $(1,0)$ to $(mq+1,q)$ passes through the lattice
points $(m+1,1)$, $(2m+1,2)$, $(3m+1,3)$, and so on, then slide the left
endpoint from $(1,0)$ to $(0,0)$). This equals the sequence obtained for
the middle tunnel using the continued fraction expansion $(mq+1)/q=[m,q]$,
and Theorem~\ref{thm:middle_tunnels_slopes} also gives all $s_i=0$. For the
lower tunnel, the sequence is $[1/3]$, $3,\ldots\,$ $3$, $5,\ldots\,$ $5$,
$7,\ldots\,$ $7\,\ldots\,$, $2q-1$, where each value is repeated $m$ times,
except that $3$ appears $m-1$ times. Thus the middle tunnel is equivalent
to the upper tunnel and distinct from the lower tunnel.

For the case when $p=mq-1$, a similar examination (using the line segment
from $(0,0)$ to $(mq,q)$ and sliding the right-hand endpoint to $(mq-1,q)$)
finds the slopes to be $[1/(2m-1)]$, $4m-1,\ldots\,$, $2m(q-1)-1$.  The
continued fraction expansion is $(mq-1)/q=[m-1,1,q-1]$, and the algorithm
for the middle tunnel gives the same slope sequence and all $s_i=0$. For
the lower tunnel, the sequence is $[1/3]$, $3,\ldots\,$ $3$, $5,\ldots\,$
$5$, $7,\ldots\,$ $7\,\ldots\,$, $2q-1$, where each value is repeated $m$
times, except that $3$ and $2q-1$ are repeated $m-1$ times. Again, the
middle tunnel is equivalent to the upper tunnel and distinct from the lower
tunnel.
\smallskip

\noindent\textsl{Case III.} Neither Case~I nor Case~II
\smallskip

In these cases, Theorem~\ref{thm:middle_tunnels_slopes} shows that the
middle tunnel has at least one nonzero value of $s_i$, so is distinct from
the upper and lower tunnels. Reducing to the case when $p>q\geq 2$,
Theorem~\ref{thm:cabling_sequence} shows that the slopes are all distinct
for the upper tunnel, but there is a repeated slope for the lower
tunnel. This completes the verification.

We note that the cases when there are fewer than three tunnels are exactly
those for which the middle tunnel is semisimple. This verifies the
equivalence of the first two conditions in the following proposition. The
equivalence of the first and third is from~\cite{B-R-Z}.

\begin{proposition} For the
$(p,q)$ torus knot $K_{p,q}$, the following are equivalent:
\begin{enumerate}
\item $p\not\equiv \pm1 \pmod q$ and $q\not\equiv \pm1 \pmod p$.
\item The middle tunnel is regular.
\item $K_{p,q}$ has exactly three tunnels.
\end{enumerate}
\label{prop:torus_regular}
\end{proposition}

\bibliographystyle{amsplain}

\end{document}